\documentclass[11pt]{llncs}

\usepackage{url}

\usepackage{amssymb}
\usepackage{enumerate}

{\bfseries}{\itshape}
{\bfseries}{\itshape}
{\bfseries}{\itshape}
{\bfseries}{\itshape}
{\bfseries}{\itshape}

\newcommand {\bbox}{\rule{0.6em}{0.6em}}

\title{On the boxicity and cubicity of hypercubes}
\author{L. Sunil Chandran 
\thanks{Indian Institute of Science,
Dept. of Computer Science and Automation,
Bangalore  560012, India.  email: \emph{sunil@csa.iisc.ernet.in}}
\and Naveen Sivadasan
\thanks{Strand Genomics,
237, Sir. C. V. Raman Avenue, Rajmahal Vilas,
Bangalore  560080, India.  email: \emph{naveen@strandgenomics.com}}
\institute{}}

\date{}

\begin{document}

\pagestyle{plain}
\pagenumbering{arabic}
\maketitle

\newcommand{\degr}{\ensuremath{\Delta}}
\newcommand{\boxi}{\ensuremath{\mathrm{box}}}
\newcommand{\sph}{\ensuremath{\mathrm{sph}}}
\newcommand{\cubi}{\ensuremath{\mathrm{cub}}}
\newcommand{\neigh}{\ensuremath{N}}
\newcommand{\dist}{\ensuremath{\delta}}
\newcommand{\ignore}[1]{}

\begin {abstract}
For a  graph $G$, its \emph{cubicity}  $\cubi(G)$ is the minimum dimension $k$ such that
$G$ is representable as the intersection graph of 
(axis--parallel) cubes in $k$--dimensional
space.
Chandran et al. \cite{CMO} showed that for a $d$--dimensional hypercube $H_d$,~
$\frac{d-1}{\log d} \le \cubi(H_d) \le 2d$.
In this paper, we show that $\cubi(H_d) = \Theta\left( \frac{d}{\log d}\right)$.
The parameter \emph{boxicity} generalizes cubicity: the boxicity $\boxi(G)$ of a graph $G$
is defined as the minimum dimension $k$ such that $G$ is representable
as the intersection graph of axis parallel boxes in $k$ dimensional space.
Since $\boxi(G) \le \cubi(G)$ for any graph $G$, our result implies that
$\boxi(H_d) = O\left(\frac{d}{\log d}\right)$. The problem of determining 
a non-trivial lower bound for $\boxi(H_d)$ is left open.
\end {abstract}

\section {Introduction}

Let ${\mathcal{F}} = \{ S_x \subseteq U : x \in V \}$ be a family of 
subsets of a universe $U$, where $V$ is an index set. The intersection graph 
$\Omega({\mathcal{F}})$ of ${\mathcal{F}}$ 
has $V$ as vertex set, and two distinct vertices
$x$ and $y$ are adjacent if and only if $S_x \cap S_y \ne \emptyset$.
Representations of graphs as the intersection graphs of various
geometrical objects is a well studied topic in graph theory.  Two well-known
concepts in this area of graph theory are the cubicity and the boxicity.
These concepts were introduced by F. S. Roberts in 1969  \cite {Roberts} and
they  find applications in niche overlap in ecology and to
problems of fleet maintenance in operations research. (See \cite {CozRob}.) 
 
A $k$--dimensional box is a Cartesian product 
$R_1 \times R_2 \times \cdots \times R_k$ where $R_i$ (for $1 \le i \le k$)
is a closed interval of the form $[a_i, b_i]$ on the real line. A 
$k$--dimensional cube is a Cartesian product $R_1 \times R_2 \times \cdots
\times R_k$, where $R_i$ is a closed interval of the form $[a_i, a_i + 1]$
on the real line. 
For a  graph $G$, its \emph{boxicity} is the minimum dimension $k$, such that
$G$ is representable as the intersection graph of 
(axis--parallel) boxes in $k$--dimensional
space. We denote the boxicity of a graph $G$ by $\boxi(G)$. 
When the boxes are restricted to  be (axis--parallel) $k$--dimensional 
cubes, the minimum dimension $k$ required 
to represent $G$ is called the \emph{cubicity}
of $G$ and is denoted by $\cubi(G)$.  
It is easy to see that for any graph $G$, 
$ \boxi (G) \le \cubi(G)$.

A $d$--dimensional hypercube $H_d$ on $2^d$ vertices is defined as follows.
The vertices of $H_d$ correspond to 
the $2^d$  binary strings  each of length $d$,
two of the vertices  
being adjacent \emph{if and only if} the corresponding binary strings differ
 in exactly one bit position.
Hypercubes are a well-studied class of graphs, which arise in the
context of parallel computing, coding theory, algebraic graph theory
and many other areas. Hypercubes are popular among graph theorists
because of their symmetry, small diameter and many other interesting
graph--theoretic properties.

It was shown by Chandran, Mannino and Oriolo   \cite{CMO} that
$\frac{d-1}{\log d} \le \cubi(H_d) \le 2d$.  
In this paper, we show the following:
$$
 \cubi(H_d) = \theta\left( \frac{d}{\log d} \right) \enspace.
$$
Since $\cubi(G)$ is an upper bound for $\boxi(G)$, clearly the above
result also implies that $\boxi(H_d) \le \frac {c d}{\log d}$ where $c$ is
a constant. Such an upper bound for $\boxi(H_d)$ was not known before. 
We leave open the question of determining a non-trivial lower bound for $\boxi(H_d)$.

\ignore{
	\subsection {Our result} 
	
	The $d$--dimensional hypercube $H_d$ is the graph on 
	$2^d$ vertices, which correspond to the $2^d$  $d$--vectors  
	whose components are either 0 or 1, two of the vertices  
	being adjacent \emph{if and only if} they differ in just one coordinate.
	Hypercubes are a well-studied class of graphs, which arise in the
	context of parallel computing, coding theory, algebraic graph theory
	and many other areas. Hypercubes are popular among graph theorists
	because of their symmetry, small diameter and many interesting
	graph--theoretic properties.
	
	Chandran, Mannino and Oriolo \cite {CMO} showed that $\frac {d}{\log d} \le 
	\cubi(H_d) \le 2d$. In this paper we show the following: 
	 $$ \cubi(H_d) = \theta \left ( \frac {d}{\log d} \right )$$
	In a light sense, this seems like a nonintuitive result, since one tends to
	expect that the cubicity of a $d$-dimensional hypercube may be $\theta(d)$. 

}

\subsection {A brief literature survey on cubicity  and boxicity}

It was shown 
by Cozzens \cite {Coz}  that computing the boxicity of a graph is NP--hard. 
This
was later improved by Yannakakis \cite {Yan1}, 
and finally by  Kratochvil \cite {Kratochvil} 
who  showed that deciding  whether
boxicity of a graph is at most 2 itself 
is NP--complete. The complexity of finding
the maximum independent set in bounded boxicity graphs was considered
by \cite {ImaiAsano,Fowler}.

There have also been attempts to estimate or bound the boxicity of 
graph classes  with special structure. 
Scheinerman \cite {Scheiner} 
showed that the boxicity of outer planar graphs is at most $2$.
Thomassen \cite {Thoma1} proved that 
the boxicity of planar graphs is
bounded above by $3$. Upper bounds for the boxicity of
many other graph classes such as chordal graphs, AT-free graphs, permutation graphs etc.,
were shown in \cite{CN05} by relating the boxicity of a graph with its treewidth.

Researchers have also tried to generalize or extend  the
concept of boxicity in various ways. The poset boxicity \cite {TroWest}, 
the rectangle number \cite {ChangWest}, grid dimension \cite {Bellantoni},
circular dimension \cite {Feinberg,Shearer}  and the boxicity 
of digraphs \cite {ChangWest1} are some  examples.

\section {Definitions and Notations}
Let $G$ be a undirected simple graph. 
We denote by $V(G)$ and $E(G)$ the vertex and edge sets 
of $G$, respectively.

As mentioned in the introduction, a  string of length $d$ consisting
 only of 0s and 1s (i.e., a {\it binary string}) 
can be associated (in one-to-one correspondence)
 with each vertex of a $d$--dimensional hypercube $H_d$, such that
two vertices $u$ and $v$ are adjacent if and only if their
corresponding binary strings differ in exactly one position.
Let $f(v)$ denote the binary string associated with  the vertex $v$.
The value of the binary digit (i.e., {\it bit} )  at the $i$--th position 
of the binary string $f(v)$ will
be denoted by $f_i(v)$. 

Given two vertices $u$ and $v$, let $D(u,v) = \{ i : 1 \le i \le d \mbox { and }
  f_i(v) \ne f_i(u) \}$. That is $D(u,v)$ is the set of ``positions'' where
 the bit values of  $f(u)$ and $f(v)$ differ from each other. The Hamming distance
between $f(u)$ and $f(v)$ is defined to be $|D(u,v)|$. 
It is easy to observe that the shortest distance $\dist(u,v)$ between two vertices
$u$ and $v$ equals the Hamming distance between $f(u)$ and $f(v)$. That is,
$\dist(u,v) = |D(u,v)|$. 

~~~~~

\noindent {\bf Unit Interval graphs:} A graph $G$ is a unit interval 
graph if and only if each
vertex of  $G$ can be mapped to a (closed) interval of unit length 
on the real 
line such that two distinct vertices are adjacent in $G$ if and only if 
the corresponding (unit) intervals intersect. 

~~~~

\noindent The following  characterization of cubicity is easy to prove. (See
\cite {Roberts}.)

~~~~

\begin {lemma}
\label {EquivalentCharacterCubicityLemma}
Let  $G$ be a simple graph. 
Let $t$ be the minimum integer such that there exists $t$ unit interval graphs
$G_1,G_2, \cdots, G_t$  on the same vertex  set as that of $G$ (i.e., $V(G_i)
=V(G)$ for each $i$, $1 \le i \le t$),
such that $E(G) = E(G_1) \cap E(G_2) \cap \cdots 
\cap E(G_t)$. Then, cubicity$(G)= t$.
\end {lemma}

\section {Upper bound for $\cubi(H_d)$}

In this section we will show that there exists a constant $c$, such that
$\cubi(H_d) \le \frac {cd}{\log d}$. 
By  Lemma \ref {EquivalentCharacterCubicityLemma},
it is sufficient to 
demonstrate that there exist $\frac {c.d}{\log d}$
unit interval graphs  (where $c$ is a constant) on the same vertex set as that of $H_d$,
 such that the edge set of $H_d$ is the intersection
of the edge sets of these unit interval graphs. With this in mind,
corresponding to each vertex $x \in V(H_d)$, 
we define below a special unit interval graph $I_x$.

~~~~~~~

\noindent {\bf Construction of the unit interval graph $I_x$:} We map
$x$ to the interval $[0,1]$. Let $u \in V(H_d) - \{x\}$. 
We map $u$ to  the interval $[\dist(x,u),\dist(x,u) +1]$. (Recall that
$\dist(x,u)$ denotes the shortest distance between $u$ and $x$ in $H_d$.)
Let $I_x$ be the resulting interval graph (with vertex set  $V(H_d)$).

\begin {lemma}
For any $x \in V(H_d)$, $E(I_x) \supseteq E(H_d)$.  
\end {lemma} 

\proof { Let $(u,v) \in E(H_d)$. Without loss of generality, assume that
   $\dist(x,v) \ge\dist(x,u)$.  Note that $H_d$ is a bipartite graph.
   Hence, if $(u,v)$ is an edge of $H_d$ then $\dist(x,u) \ne\dist(x,v)$. 
   (Otherwise, there will be an odd cycle in $H_d$.)  Moreover,
   $\dist(x,v) \le\dist(x,u) + 1$, since $(u,v) \in E(H_d)$. It follows that 
   $\dist(x,v) =\dist(x,u) + 1$.  Thus the intervals associated with $u$ and
   $v$ in $I_x$ intersect. (They touch each other at $\dist(x,v)=\dist(x,u) + 1$.)
   Hence the Lemma.  
   }$\bbox$

Our plan is to show that there exists a subset $S \subset V(H_d)$
 with $|S| \le  \frac {c d}{\log d}$ (where $c$ is a constant) such 
that  $\bigcap_{x \in S} E(I_x) = E(H_d)$. The reader may note that
in view of Lemma \ref {EquivalentCharacterCubicityLemma},
it is sufficient to show that such a set $S$ has the following
property. (We name this property as Property $P$.) 

\begin {definition}[Property P]
A subset $S \subset V(H_d)$ is said to have the property
$P$ if and only if  for each  $(u,v) \notin E(H_d)$,  there exists 
a vertex $x \in S$, such that $(u,v) \notin E(I_x)$. 
\end {definition} 

\noindent {\bf Choosing the subset $S$ randomly:} 
We select a random subset $S$ of $V(H_d)$ by conducting the 
following experiment:  We select a binary string $x$ such that
the bit at position $i$ is set to $1$ with probability $\frac {1}{2}$. That
is, for any $i$, $1 \le i \le d$, $Pr(f_i(x) = 0 ) = \frac {1}{2}$ and
$Pr(f_i(x) = 1) = \frac {1}{2}$. We do this experiment $\frac {cd}{\log d}$
times, thus selecting $\frac {cd}{\log d}$ binary strings. Let $S$ be the 
multi--set of vertices which correspond to the strings so selected. 
Clearly, $|S| = \frac {cd}{\log d}$. 

~~~~~~~~

We show 
that if subset $S$ is constructed randomly as explained above,
then $Pr( S  \mbox { doesn't satisfy property P } ) < 1$. 
As a consequence, it follows that there exists a subset 
$S$ of $V(H_d)$, where $|S| \le \frac {cd}{\log d}$ ($c$ being a constant),
such that $S$ satisfies property $P$.

~~~~~~~~

\noindent The following Lemma is an easy consequence of the construction of $I_x$. 
 
\begin {lemma}
\label {EdgeDistanceRelationLemma} 
For any vertex $x \in V(H_d)$, $(u,v) \in E(I_x)$ if and only if 
$|\dist(x,u) -\dist(x,v)| \le  1$.
\end {lemma} 

Given three vertices  $x,u,v \in V(H_d)$, we partition the bits of 
$f(x)$ in to three categories:  

\begin {enumerate} 

\item  If $i \notin D(u,v)$, then $f_i(x)$ is defined to be  a 
 \emph{neutral} bit of $f(x)$ (with respect to $u$ and $v$). ({\it The
 reason why we name $f_i(x)$ a  neutral bit is the following: 
 If $i \notin D(u,v)$, then it is the case that either 
 $i \in D(u,x)$ and $i \in D(v,x)$ or   
 $i \notin D(u,x)$ and  $i \notin D(v,x)$. })

\item If $i \in D(u,v)$, then $f_i(x)$ is called a $u$--bit if and only 
 if $f_i(x) \ne f_i(u)$.  Clearly in that case $f_i(x) = f_i(v)$.
 {\it Reader may note  that if $f_i(x)$ is a $u$--bit  then $i \in D(x,u)$,
whereas $i \notin D(x, v)$.}
The number of $u$--bits of $f(x)$ will be denoted by $n_u(x)$.  

\item If $i \in D(u,v)$, then $f_i(x)$ is called a $v$--bit if and only
 if $f_i(x) \ne f_i(v)$. Clearly if $i \in D(u,v)$, $f_i(x)$ is a 
 $v$--bit if and only if it is not a $u$--bit. {\it It may be noted
 that if $f_i(x)$ is a $v$--bit then $i \in D(x, v)$, whereas $i \notin D(x, u)$. }
 The number of $v$--bits of $f(x)$ will be denoted by $n_v(x)$. 

\end {enumerate}

\noindent The next lemma follows immediately, from the discussion above. 

\begin {lemma}
\label {DistanceBitRelationLemma} 
Let $x,u,v \in V(H_d)$. Then $|D(u, v)| = n_u(x) + n_v(x)$ and
$|\dist(x,u) -\dist(x,v)| = |n_u(x) - n_v(x)|$. 
\end {lemma}

Let $x$ be a vertex corresponding to a randomly chosen binary string:
i.e., $Pr(f_i(x) = 1) = \frac {1}{2}$. 
We now bound for a pair of nonadjacent vertices $(u, v)$, $Pr( (u,v) \in E(I_x) )$ as follows.
By Lemma \ref 
{EdgeDistanceRelationLemma} and Lemma \ref {DistanceBitRelationLemma}, for
a pair of nonadjacent vertices $(u,v)$, $Pr( (u,v) \in E(I_x) ) =
 Pr( |n_u(x) - n_v(x)| \le 1)$. 
We consider two cases: 


~~~~~~~~~~~~~~~~~~~~~

\noindent
\emph{Case 1:
$\dist(u,v) = r$ is even. }
Since $u$ and $v$ are nonadjacent, $r \ge 2$. 
Since $r=\dist(u,v) = |D(u,v)|=  n_u(x) + n_v(x)$ by Lemma \ref{DistanceBitRelationLemma},
clearly $n_u(x) - n_v(x)$ is  also even.   Thus we have 
$Pr( (u,v) \in E(I_x) ) = Pr( |n_u(x) - n_v(x)| = 0)$. 
Noting that for any 
$i \in D(u,v)$, $f_i(x)$ is a $u$--bit with probability $\frac {1}{2}$ 
and it is a $v$--bit with probability $\frac {1}{2}$, we have:  

\begin {eqnarray}
\label {EvenRInequality}
Pr( (u,v) \in E(I_x) ) =  {r \choose {r/2} } 2^{-r}.
\end  {eqnarray} 

\noindent
Since $r$ is even and $r \ge 2$,  ${r \choose \frac {r}{2} } 2^{-r}
 > {{r +2} \choose \frac {r+2} {2} } 2^{-(r+2)}$. It follows that
$Pr((u,v) \in E(I_x))$ is maximized at $r =2$ and thus
$Pr((u,v) \in E(I_x)) \le \frac{1}{2}$.


~~~~~~~~~~~~~~~~~

\noindent
\emph{Case 2: $\dist(u,v)=r$ is odd. }
Since $u$ and $v$ are nonadjacent, $r \ge 3$. 
Clearly, $n_u(x) - n_v(x)$ is  odd. 
Thus we have 
$Pr( (u,v) \in E(I_x) ) = Pr (|n_u(x) - n_v(x)| = 1)$.  It follows that:

\begin {eqnarray}
\label {OddRInequality}
Pr( (u,v) \in E(I_x) )  &=& Pr \Big( n_u(x) - n_v(x) = 1\Big) + Pr\Big(n_u(x) - n_v(x) = -1\Big) \nonumber \\&=& \left ({r \choose {(r+1)/2} } + {r \choose {(r-1)/2} } \right ) 2^{-r} \nonumber \\ 
 &=&   {r \choose {(r+1)/2} } 2^{-(r-1)} 
\end {eqnarray} 

\noindent
Since $r$ is odd and $r \ge 3$,
${r \choose \frac {r+1}{2} } 2^{-(r-1)} > {{r+2} \choose {\frac {r+3}{2}}} 2^{-(r+1)}$.  
It follows that  $Pr ( (u,v) \in E(I_x) )$ is maximized at $r=3$ and thus
$Pr ( (u,v) \in E(I_x) ) \le \frac{3}{4}$.

From the above two cases, it follows that

\begin {eqnarray}
\label {SimplifiedProbabilityValueInequality} 
Pr( (u,v) \in E(I_x))  \le  \frac {3}{4}
\end {eqnarray}

\noindent
Since each $x \in S$ is chosen independently and uniformly at random,

\begin {eqnarray}
\label {WorstCaseProbInequality}
Pr( (u,v) \in E(I_x), \forall x \in S)   \le \left ( \frac {3}{4} \right )^ {|S| }  \le \left ( \frac {3}{4} \right )^ {\frac {cd}{\log d} } 
\end {eqnarray}

The obvious next step in order to derive an upper bound 
for Pr($S$ does not satisfy property $P$) would be to use the union bound, that is,
 Pr($S$ does not satisfy property $P$) $\le \sum_{(u,v) \notin E(H_d)} Pr( (u,v) \in E(I_x), \forall x \in S) $.
Unfortunately, there are ${2^d \choose 2} - 2^{d-1} d = O(2^{2d})$ 
nonadjacent pairs of vertices in $H_d$, and a straightforward application 
of the union bound as above would not suffice: the bound given by Inequality
\ref {WorstCaseProbInequality} is too weak.   But, by examining the Inequalities 
\ref {EvenRInequality} and \ref {OddRInequality}  more carefully, 
the reader can easily see that as $r$ becomes larger,
the probability that a nonadjacent pair
$(u,v)$ with $\dist(u,v) = r$ being adjacent in $I_x$ reduces, and for sufficiently
large $r$, this probability can be much smaller than what is guaranteed
by Inequality \ref {SimplifiedProbabilityValueInequality}. 
In fact, by applying 
Sterling's approximation (i.e., $n! \sim (n/e)^n \sqrt{2\pi n} $) on  Inequalities \ref {EvenRInequality} and
\ref {OddRInequality}, it is easy to verify that, there exists a constant
$c_1$, such that for a pair of nonadjacent vertices $(u, v)$, 

\begin {eqnarray}
\label {StregthenedProbabilityInequality}
Pr ( (u,v) \in E(I_x) ) \le  \frac {c_1}{\sqrt{\dist(u, v)}} 
\end {eqnarray} 

Based on this observation,
we partition  the nonadjacent pairs of vertices in $H_d$ into two  groups
$A$ and $B$  as follows:  

$$ A = \{ (u,v) : u,v \mbox { are nonadjacent in $H_d$ and }\dist(u,v) > \frac {d}{\log^2 d } \} $$

$$ B = \{ (u,v): u,v \mbox { are nonadjacent in $H_d$ and }\dist(u,v) \le \frac {d}{\log^2 d} \} $$

\begin {definition} 
 A subset $S$ of $V(H_d)$ is said to satisfy Property $P_A$ 
 (respectively $P_B$) if and only if for each nonadjacent pair
 $(u,v) \in A$ (respectively in $B$), there exists a vertex $x \in S$, 
 such that $(u,v) \notin E(I_x)$. 
\end {definition} 

\noindent It is easy to see the following: 
\begin {eqnarray}
\label {ProbabilitySumInequality} 
\hspace*{-1cm} Pr \mbox {(S does not satisfy P)} & \ \ \ \ \le \ \ \  & 
  Pr (\mbox {S does not satisfy}  P_A)   \nonumber  \\ & \ \ \   &   \ \ \ \ \ \ \  + \ Pr (S \mbox { does not satisfy  } P_B)  
\end {eqnarray}
We will show that each of the two terms in the right hand side is strictly
less than $\frac {1}{2}$, so that the left hand side is strictly less than
$1$, as required.  

~~~~~~~~~~~~~~~~~~~

\noindent
Since $|A| \le 2^{2d}$, and recalling that 
for any pair $(u,v) \in A$, we have $\dist(u, v) > \frac{d}{\log^2 d}$, we can apply union bound to
show that, 

\begin {eqnarray}
Pr (S \mbox { does not satisfy } P_A) & \le & \sum_{(u, v) \in A} \left( \frac{c_1}{\sqrt{\dist(u,v)}} \right)^{\frac{cd}{\log d}} \nonumber \\
&\le& \left ( \frac {c_1 \log d}{\sqrt{d}} \right )^{\frac {cd}{\log d} } 2^{2d} \nonumber \\
     &\le&   2^{2d + \left(\log c_1 + \log\!\log d - \frac{\log d}{2} \right)\frac{cd}{\log d}} \nonumber \\ 
   &\le&   2^{-\frac{c}{8}d} ~<~ \frac{1}{2}\enspace,
\end {eqnarray}  
when $c$ is a suitably large constant and when $d \ge c_3$ for a sufficiently large constant $c_3$.
(For a sufficiently large constant $c_3$ with $d \ge c_3$, $\frac{\log c_1 + \log\!\log d}{\log d} \le \frac{1}{4}$.
Also, for a suitably large constant $c$, $2d \le \frac{cd}{8}$.)


~~~~~~~~~~~~~~~~~~~~~

Now we deal with the pairs in $B$.
Recall that an  upper bound for $Pr ( (u,v) \in E(I_x)
\forall x \in S)$ is given  by Inequality \ref {WorstCaseProbInequality}.
But, unfortunately $|B|$ is too  big  to infer that
$Pr (S \mbox { does not satisfy $P_B$}) < \frac {1}{2},$ by a simple
application of union bound. 
To overcome this difficulty, we define an equivalence relation
$\mathcal{R}$ on $B$ such that the pairs in the same equivalence
class behaves identically,  i.e. if $(u_1,v_1)$ and $(u_2,v_2)$ belong
to the same equivalence class then for any $x \in V(H_d)$, $(u_1,v_1) 
 \in E(I_x)$ if and only if $(u_2,v_2) \in E(I_x)$. 

 Recall that $f(u)$ denotes the binary string associated with $u$.  Let
 ${\cal P} = \{k_1, k_2, \cdots, k_i\}$ where 
 $1 \le k_1 < k_2 < \cdots < k_i \le d$.  We denote by $f_{{\cal P}}(u)$
 the binary string obtained by concatenating the bits $f_{k_1}(u),f_{k_2}(u), \cdots, f_{k_i}(u)$ in that order. We call $f_{\cal P}(u)$ as the \emph{bit pattern} of $u$ at the set of positions ${\cal P}$.

From now on, for any pair  of vertices $u$ and $v$, 
we choose to represent it by the 
ordered pair $(u, v)$ if  $f_{D(u,v)}(u)$ is less than $f_{D(u,v)}(v)$ in the lexicographic order;
else we represent it by $(v, u)$. (The reader may observe that the 
bit pattern  $f_{D(u,v)}(u)$ is the complement of 
the bit pattern $f_{D(u,v)}(v)$.)


We define the equivalence relation $\mathcal{R}$ as follows:
Consider two pairs $(u_1, v_1)$ and $(u_2, v_2)$.
\begin{eqnarray*}
(u_1, v_1) \, \mathcal{R} \, (u_2, v_2)  ~~ \Leftrightarrow ~~ D(u_1, v_1) = D(u_2, v_2) \mbox{~~~~~and~~} \\ 
f_{D(u_1,v_1)}(u_1) = f_{D(u_1,v_1)}(u_2) 
\end{eqnarray*}
That is, $(u_1, v_1)$ and $(u_2, v_2)$ are related by $\mathcal{R}$ if and only if:
\emph{1)} the set of bit positions where $u_1$ differs from $v_1$ is identical to the set of positions 
where $u_2$ differs from $v_2$ and \emph{2)} the bit pattern of  $u_1$ and $u_2$ at those bit positions are
identical.

Let $B_1, \ldots, B_\alpha$ be the equivalence classes of $B$ under $\mathcal{R}$.
Note that each equivalence class $B_k$ corresponds to a unique pair $({\cal P}, s)$, where ${\cal P}$
is a set of $i$ distinct bit positions, where $2 \le i \le \frac{d}{\log^2 d}~,$ and $s$ is a binary string of length $i$. It is easy to see that the number of equivalence classes $\alpha$ has the following upper bound.
Let $t =  \left \lfloor \frac {d}{\log^2 d}\right \rfloor$. Then,
\begin{eqnarray} \label{eq:alpha}
\alpha ~\le~ \sum_{i=2}^{i=t} {d \choose i} 2^i   ~\le~ t {d \choose t}
         2^t  ~\le~ t (2d)^t 
\end{eqnarray}

Now, from the definition of the relation $\mathcal {R}$,
it is easy to see that if $(u_1, v_1)$ and $(u_2, v_2)$ are in the same
equivalence class $B_k$ then 
for any $x \in V(H_d)$, $|n_{u_1}(x) - n_{v_1}(x)|= |n_{u_2}(x) - n_{v_2}(x)|$
and therefore $(u_1,v_1) \in E(I_x)$ if and only if $(u_2, v_2) \in E(I_x)$.

\noindent Thus applying the union bound using (\ref{WorstCaseProbInequality})
 and using the inequality (\ref{eq:alpha}) for $\alpha$,
 we get:

\begin{eqnarray}
 &Pr (S \mbox { does not satisfy } P_B) ~~~~~~~~~~~~~~~~~~~~~~~~~~~~~~~~~~~~~~~~~~~~~~~~~~~~~~~~~~~~~~~~   \nonumber \\
 & =  Pr \Big( \exists (u,v) \in B : \mbox { such that }  (u,v) \in E(I_x), \forall
 x \in S \Big) \nonumber \\ &\le~  \alpha \left(\frac{3}{4}\right)^{\frac {cd}{\log d} }   ~\le~  t (2d)^t \left(\frac{3}{4}\right)^{\frac {cd}{\log d} }   ~<~  \frac {1}{2} ~~~~~~~~~~~~~~~~~~~~~~~~~~~~~~~~~~~~~~~~
\end{eqnarray}
for a suitably large constant $c$.

\ignore{

	It follows that inequality \ref{WorstCaseProbInequality} can be further strengthened to obtain the following:

	~~~~~~~~~~~

	\noindent
	For any fixed subset $X$ of $\{1, \ldots, d\}$ where $|X| > 1$,

	\begin {eqnarray}
	\label {WorstCaseProbInequality2}
	Pr\Big( \exists (u,v) : D(u, v) = X \mbox{~and~} (u, v)  \in E(I_x), \forall x \in S\Big)   \le \left ( \frac {3}{4} \right )^ {|S| }  
	\end {eqnarray}

	\noindent
	Hence 
	\begin{eqnarray*}
	&&Pr\Big( \exists (u, v) : |D(u,v)| > 1 \mbox{~and~} (u,v) \in E(I_x), \forall x \in S\Big)    \\
	& = & \sum_{X \subseteq \{1, \ldots, d\}, |X| > 1} 
	Pr\Big( \exists (u,v) : D(u, v) = X \mbox{~and~} (u, v)  \in E(I_x), \forall x \in S\Big)   \\
	& \le & \sum_{X \subseteq \{1, \ldots, d\}, |X| > 1} \left(\frac {3}{4} \right )^ {|S| }  
	\end{eqnarray*}

	to get an upper bound for $Pr(S \mbox {does not satisfy property} P_B)$
	we need to consider the  contribution from only one such pair. Motivated
	by this observation, we partition the set $B$ into a collection of subsets
	so that all the pairs within a given subset behave {\it identically}: if
	one of them belongs to $E(I_x)$, then all others also belong to $E(I_x)$
	and vice versa. We partition $B$ step by step as given below:  

	 \begin {enumerate} 
	 \item 
	 For  $i = 2, \cdots,\left \lfloor  \frac {d}{\log^d} \right \rfloor $, 
	 let $B_i \subseteq B$ be defined as: 
	 $$B_i = \{ (u,v) :\dist(u,v) = i\}$$
	Clearly $\bigcup_{i=2}^{i = \left \lfloor  \frac {d}{\log^d} \right \rfloor} =B$

	 \item  We partition  each $B_i$ further as follows:
	 Let $P_i$ represent the set of  all combinations of $i$ numbers selected from
	 $\{1,2, \cdots, d\}$. Clearly  $|P_i| = {d \choose i}$ and any element
	 $p \in P_i$ is a set of the form $\{k_1,k_2, \cdots, k_i \}$ where 
	 $1 \le k_j \le d$, for $1 \le i \le i$. 
	 For each  $p \in P_i$, we define 
	 $$B_{i,p} = \{ (u,v) \in B_i : D(u,v) = p\}$$ 
	 Noting that any pair $(u,v) \in B_i$ is such that $f(u)$ and $f(v)$ differs
	 in exactly $i$ positions (by the definition of $B_i$), clearly the 
	 collection of subsets $B_{i,p}$ partition $B_i$. 
	 
	 \item Recall that $f(u)$ denotes the binary string associated with $u$.  Let
	 $p \in P_i$, and let $p = \{k_1, k_2, \cdots, k_i\}$ where 
	 $1 \le k_1 < k_2 < \cdots < k_i \le d$.  We denote by $f_{[p]}(u)$
	 the binary string obtained by concatenating the bits $f_{k_1}(x),f_{k_2}(x), \cdots, f_{k_i}(x)$ in that order.  Let $S_i$ represent
	 the set of all  distinct binary strings of length $i$.
	 For a binary string $s$, let $\overline s$ denote its complement. 
	 Now for  each $s \in S_i$ we define: 
	 $$ B_{i,p,s} = \{ (u,v) \in B_{i,p} : f_{[p]}(u) = s \mbox { or } \overline s
	 \}$$ 

	\end {enumerate} 
	 
	 In the following lemma we show that all the pairs in any $B(i,p,s)$
	 behaves {\it identically}  with respect to any fixed vertex  $x \in V(H_d)$.
	 That is either all the pairs in $B(i,p,s)$ are edges of $I_x$ or 
	 none of them are  edges of $I_x$.  

	 \begin {lemma}
	 For any $x$,  either $B(i,p,s) \subseteq E(I_x)$ or $B(i,p,s) \cap
	 E(I_x) = \emptyset$.
	 \end {lemma} 

	 Finally we note that the number of sets of the form $B(i,p,s)$,
	 where $1 \le i \le x$ where 
	 $x =  \left \lfloor \frac {d}{\log^2 d}\right \rfloor$, is 
	 $$ \sum_{i=2}^{i=x} {d \choose i} 2^i   \le x {d \choose x}
		 2^x  \le x (2d)^x $$ 

	\noindent Thus applying the union bound (and choosing $c$ to be suitably large)  we get:

	\begin {eqnarray}
	 Pr ( \exists (u,v) \in B : \mbox { such that }  (u,v) \in E(I_x), \forall
	 x \in S ) \le  2^{\frac {-cd}{\log d} } x (2d)^x  \le \frac {1}{2}
	\end {eqnarray}

}

~~~~~~~~~~~~~~

Thus recalling inequality \ref {ProbabilitySumInequality}, we
have $Pr$($S$ does not satisfy property $P$) $ < 1$. It follows that there
exists a subset $S \subseteq V(H_d)$, with $|S| \le \frac {cd}{\log d}$,
such that $S$ satisfies property $P$. In other words: 

\begin {theorem} 
\label {UpperboundTheorem}
$\cubi(H_d) \le  \frac {cd}{\log d}$, where $c$ is a
constant.
\end {theorem}

\noindent The following lower bound for the cubicity of $H_d$ 
was shown in \cite {CMO}. 

\begin {theorem} [Chandran et al. \cite{CMO}]
\label {LowerboundThm}
$\cubi(H_d) ~\ge~ \frac {d-1}{\log d} \enspace.$
\end {theorem}

Finally combining the  upper bound of Theorem 
\ref {UpperboundTheorem}  with the lower bound of Theorem \ref {LowerboundThm},
we have:

\begin {theorem}
$\cubi(H_d) = \theta \left ( \frac {d}{\log d} \right ).$
\end {theorem}


\end {document}